\documentclass[10pt]{amsart}

\usepackage[dvips]{epsfig}
\usepackage{amscd}
\usepackage{amssymb}
\usepackage{amsthm}
\usepackage{amsmath}
\usepackage{latexsym}
\usepackage[all]{xy} 

\theoremstyle{plain}
\newtheorem{theorem}{Theorem}[section]
\theoremstyle{plain}

\newtheorem{lemma}[theorem]{Lemma}
\newtheorem{proposition}[theorem]{Proposition}
\newtheorem{corollary}[theorem]{Corollary}
\theoremstyle{definition}

\theoremstyle{remark}
\newtheorem{remark}{Remark}[section]

\newcommand{\A}{\ensuremath{\mathcal{A}}}
\newcommand{\B}{\ensuremath{\mathcal{B}}}
\newcommand{\eps}{\ensuremath{\varepsilon}}
\newcommand{\R}{\ensuremath{\mathbb{R}}}
\newcommand{\N}{\ensuremath{\mathbb{N}}}
\newcommand{\Z}{\ensuremath{\mathbb{Z}}}
\newcommand{\BS}{\ensuremath{\mathbb{S}}}

\newcommand{\del}{\partial}
\newcommand{\ol}{\overline}
\newcommand{\ra}{\rightarrow}
\newcommand{\alp}{\alpha}

\newcommand{\tr}{\mathrm{tr}}
\newcommand{\bx}{\bar{x}}
\newcommand{\by}{\bar{y}}
\newcommand{\bu}{\bar{u}}
\newcommand{\ba}{\bar{a}}
\newcommand{\bb}{\bar{b}}
\newcommand{\bc}{\bar{c}}
\newcommand{\bff}{\bar{f}}
\newcommand{\bs}{\bar{\sigma}}

\numberwithin{equation}{section}
\allowdisplaybreaks

\title[Error bounds]{
On error bounds for monotone approximation schemes for multi-dimensional
Isaacs equations.
} 

\date{\today}

\author[Jakobsen]{Espen R.~Jakobsen}
\address[Espen R.~Jakobsen]{\newline
    Department of Mathematical Sciences\newline
    Norwegian University of Science and Technology\newline
    N--7491 Trondheim, Norway}
\email[]{erj\@@math.ntnu.no}
\urladdr{http://www.math.ntnu.no/\~{}erj/}


\begin{document}

\begin{abstract}
Recently, Krylov, Barles, and Jakobsen developed the theory for
estimating errors of monotone approximation schemes for the Bellman
equation (a convex Isaacs equation). In this 
paper we consider an extension of this theory to 
a class of non-convex multidimensional Isaacs equations. This is the
first result of this kind for non-convex multidimensional fully
non-linear problems.  

To get the error bound, a key intermediate step is to introduce a
penalization approximation.
We conclude by (i) providing new error bounds for penalization approximations
extending earlier results by e.g. Benssousan and Lions, and (ii) obtaining
error bounds for approximation schemes for the penalization equation
using very precise a priori bounds and a slight generalization of
the recent theory of Krylov, Barles, and Jakobsen. 
\end{abstract}

\maketitle

\section{Introduction}
\label{Sec:intro}
In this paper we will study error bounds for approximation schemes for
a class of non-convex multidimensional Isaacs equations. To be
precise, 
we will consider the following (one-sided) obstacle problem
\begin{align}
\label{E}
\min\{F(x,u,Du,D^2u),u-g\}=0 \quad \text{in}\quad \R^N,
\end{align}
where $g$ is the obstacle ($u\geq g$), and $F$ is given by
\begin{align}
\label{defF}
F(x,r,p,X)=\sup_{\alp\in\A}\left\{-\tr[a^\alp(x)X]+b^\alp(x)p
  + f^{\alp}(x,r)\right\}. 
\end{align}
Here $\A$ is a compact metric space, $a$ is a positive semidefinite
matrix, $f$ is strictly increasing in $r$, and the data is a least
bounded and uniformly continuous. Precise assumptions will be
specified later.  
This equation is non-convex because of the $\min/\sup$ form, but also
because the $f$ term may be non-convex in $r$. It may also be degenerate
since $a$ may vanish at certain $x$ and $\alp$. Under the assumptions
we will use, this equations can always be rewritten as an Isaacs
equation,
\begin{align*}
\inf_{\alp\in\A}\sup_{\beta\in\B}\left\{-\tr\left[\bar a^{\alp,\beta}(x)D^2u\right]+
\bar b^{\alp,\beta}(x)Du+\bar c^{\alp,\beta}(x)u+\bar
f^{\alp,\beta}(x)\right\}=0 
\end{align*}
in $\R^N$, for suitably defined $\bar a$, $\bar b$, $\bar c$, and
$\bar f$. It is well-known that such equations do not in general have
smooth solutions.

The above problem \eqref{E} is also called a variational inequality, and such
problems occur in many applications and have been studied intensively
for a long time. The classical theory  for variational inequalities
(see e.g. \cite{KS:Book,BL:Book}) studies weak or variational solutions
and uses either PDE techniques and Sobolev space theory or
probabilistic techniques using the connection with optimal stopping
time problems. In this paper we will (mostly) consider viscosity
solutions which is a weaker and more recent
concept of solutions. 
We refer to \cite{CIL:UG,FS:Book} for a general overview of the
viscosity solution theory, and to
e.g. \cite{Ph:OS,Am:OP,BKR:BSAmOp}
for analysis and applications of obstacle problems in the viscosity
solutions setting. We mention in particular the many applications in
finance, like e.g. the pricing problem for American options \cite{WHD:Book}.

In the viscosity solution setting the first results on error bounds
for monotone schemes were obtained by Crandall and Lions
\cite{CL:Appr} for first order equations. This case has later been
studied by many authors. 
Only recently did Krylov \cite{Kr:HJB1,Kr:HJB2} obtain
the first results for second order fully non-linear equations
(the convex Bellman-equations), and these results were then extended
and improved by Barles and Jakobsen \cite{BJ:Rate,BJ:Rate2,J:Par}. We
refer to the recent paper \cite{BJ:Rate2} for the best results
available at the present time. All these results concern the convex
Bellman equation. In the non-convex fully non-linear case, there are
to the best of the author's knowledge no results in the
multi-dimensional case. The only non-convex result we know about
applies to one dimensional problems \cite{J:Non-conv}.

In this paper we will give error bounds for general monotone
approximation schemes for the non-convex multi-dimensional problem
\eqref{E}. We will use the following abstract notation for such schemes,
\begin{align}
\label{S}
\min\left\{S(h,x,u_h(x),[u_h]_x); u_h(x)-g(x)\right\}=0\quad
\text{in}\quad\R^N, 
\end{align}
where $S$ is loosely speaking a consistent, monotone, and
uniformly continuous approximation of $F$ in \eqref{E}. The
approximate solution is $u_h$, $[u_h]_x$ is a function defined from
$u_h$, and the approximation parameter is $h$. This notation
was introduced by Barles and Souganidis \cite{BS:Conv} to display
clearly the monotonicity of the scheme: $S$ is non-decreasing in $u_h$ and
non-increasing in $[u_h]_x$. Typical approximations $S$ that we have
in mind are certain finite difference methods (FDMs) \cite{KD:Book} and
so-called control schemes \cite{CF:Appr}. In Section \ref{Sec:App} we will
explain the notation for a concrete FDM. 

To get an idea of our results, we will now consider an explicit 1D problem:
\begin{align*}
&\min\left\{\sup_{\alp\in\A}\Big\{-a^\alp\Delta^2_hu''
+f^\alp(x,u)\Big\};u-g(x)\right\}=0 \quad \text{in}\quad \R^1.
\end{align*}
We approximate this problem using a monotone FDM,
\begin{align*}
&\min\left\{\sup_{\alp\in\A}\Big\{-a^\alp\Delta^2_h u_h
+f^\alp(x,u_h)\Big\};u_h-g(x)\right\}=0 \quad \text{in}\quad \R^1,
\end{align*}
where
\begin{align*}
\Delta^2_h \phi(x)&=\frac{\phi(x+h)-2\phi(x)+\phi(x-h)}{h^2}.
\end{align*}
Under suitable assumptions on $f$ and $g$ our results (cf. Proposition
\ref{LinSRate}) give 
the following error bound,
$$\|u - u_h\|_{L^\infty}\leq C h^{1/6}.$$
If we take slightly stronger assumptions on the obstacle $g$, we get
$$\|u - u_h\|_{L^\infty}\leq C h^{1/4}.$$ 
The same results hold for problems in arbitrary space dimensions, see
Section \ref{Sec:App}. In the convex case ($g\equiv-\infty$) the
corresponding rate is $1/2$ \cite{J:Par}. 

In the (much) more difficult case when $a$ also depend on $x$, the
rate in the convex case is a least $1/5$ \cite{BJ:Rate2}. In this
paper we will not be able to handle FDMs when the coefficient $a$
depends on $x$. Results in this direction can be obtained by combining
the methods of this paper with those of \cite{BJ:Rate2}. But since the
arguments become much more involved in this case, we have chosen to
omit it. However note that the results in this paper applies to the
so-called control schemes in the general case (when $a$ depends also
on $x$), see Section \ref{Sec:App}. For a better discussion of this
point we refer to \cite{BJ:Rate}.

Let us now try to explain how we get our error bounds. As a key
intermediate step we introduce the following penalization problem,
\begin{align}
\label{PE}
F(x,v_\eps,Dv_\eps,D^2v_\eps)=\frac{1}{\eps}(v_\eps-g)^-\quad \text{in}
\quad\R^N, 
\end{align}
where $(\cdot)^-=-\min(\cdot,0)$. 
Under suitable assumptions on the data, it is possible to show that
the solution of the penalized problem \eqref{PE} converges
monotonically to the solution of the obstacle problem
\eqref{E} as $\eps\ra0$ \cite{Am:OP}. In this paper we prove new error
bounds for this convergence using easy comparison arguments.

The next step is then to consider the approximation scheme associated
to \eqref{PE} via \eqref{S},
\begin{align}
\label{PS} 
S\big(h,x,v_{h,\eps}(x),[v_{h,\eps}]_{x}\big)=\frac{1}{\eps}\left(v_{h,\eps}-
g\right)^- \quad\text{in}\quad \R^N.
\end{align}
Again we prove that $v_{h,\eps}$ converge to $u_h$ the solution of
\eqref{S} with a given error bound. This argument is completely
similar to the one mentioned above in connection with \eqref{PE}.

The third and more difficult step is to obtain error bounds for
convergence the solution $v_{h,\eps}$ of \eqref{PS} to the solution
$v_\eps$ of \eqref{PE}. To get this result we use a slight extension
of the arguments in \cite{BJ:Rate,J:Par}. What is new here is that the
equation need not be convex in the zero-order term (the $u$ term), as
is the case for \eqref{PE}. 

As a finial step we combine the previous steps to get the full error
bound via the triangle inequality,
$$\|u-u_h\|_{L^{\infty}}\leq
\|u-v_\eps\|_{L^{\infty}}+\|v_\eps-v_{h,\eps}\|_{L^{\infty}}
+\|v_{h,\eps}-u_h\|_{L^{\infty}}.$$ 
The right hand side will depend on $h$ and $\eps$, and the result
follows after a minimization over $\eps>0$. Warning! This last step is {\em
  only} possible to perform if the bound on
$\|v_\eps-v_{h,\eps}\|_{L^{\infty}}$ does not depend on $\eps$ in a
too singular manner. Note that some coefficients in \eqref{PE} and
\eqref{PS} depend on $1/\eps$, and that naive computations would lead to a
priori bounds on the solutions that 
also depend on $1/\eps$. With such bounds, we would not be able to prove
any error bounds. For our purpose, we need and prove more
precise a priori bounds than can be found in the literature.

Let us now return briefly to the penalization problem
\eqref{PE}. Usually it is easier to obtain existence of 
solutions of the penalized problem and than of the corresponding
obstacle problem. The limit procedure ($\eps\ra0$) then gives existence
of also for the obstacle problem. We refer to Bensoussan and Lions
\cite{BL:Book} for the classical theory and to Amadori \cite{Am:OP}
for a viscosity solutions approach. Error bounds exist in the classical
case. E.g. in \cite[p. 197]{BL:Book} the following bound is proved,
\begin{align}
\label{ClRes}
\|v_\eps-u\|_{W^{1,2}}\leq C \eps^{1/2},
\end{align}
in the case when $F$ in \eqref{E} is linear, uniformly elliptic, and
in divergence form.
In this paper we prove under suitable assumptions that 
\begin{align}
\label{NewRes}
\|v_\eps-u\|_{L^\infty}\leq C \eps^{1/2},
\end{align}
and under sightly stronger assumptions on $g$ we get
$\|v_\eps-u\|_{L^\infty}\leq C \eps$. These results applies to very
general equations, see Section \ref{Sec:PM}, to all kinds of weak
solutions as long as the 
comparison principle holds, and even to monotone schemes like
\eqref{S} and \eqref{PS}. To the best of the author's knowledge this
result is new, even in the linear uniformly elliptic case e.g. under
the assumptions leading to \eqref{ClRes}. Also note that \eqref{NewRes}
does {\em not} follow from \eqref{ClRes} except in one space dimension
(by Sobolev embedding).

Let us now introduce some notation: We will use the following (semi) norms,
$$|f|_0=\mathrm{ess\, sup}_{x\in\R^N}|f(x)|,\
[f]_\mu=\mathrm{ess\, sup}_{x,y\in\R^N}\frac{|f(x)-f(y)|}{|x-y|^\mu},\
|f|_\mu=|f|_0+[f]_\mu,$$ 
where $f:\R^N\ra\R^M$ is a function and $\mu\in(0,1]$. The same notation
will be used for vector and matrix valued functions $f$, in which case
$|f|$ is interpreted as a vector and matrix norm respectively.
$L^\infty(\R^N)$, $C(\R^N)$, $C_b(\R^N)$, $C^{0,\mu}(\R^N)$, $\mu\in(0,1]$, 
$C^k(\R^N)$, $k\in\N$, denote the spaces of functions $f:\R^N\ra\R$
that are bounded, continuous, bounded and continuous, have finite norm
$|f|_\mu$, and are $k$-times continuous differentiable
respectively. Furthermore, $W^{1,2}$,
$W^{1,2}_0$, $W^{2,2}$, and $W^{1,\infty}=C^{0,1}$ are standard
Sobolev spaces.
The space of real symmetric $N\times N$ matrices are denoted by $\BS^N$, and
$X\geq Y$ in $\BS^N$ will mean that $X-Y$ is positive
semi-definite. Finally, by $D^k\phi$ we mean the vector of $k$-order
partial derivatives of a function $\phi$.

The outline of the rest of this paper is as follows: In the next
section, we treat the penalization problem \eqref{PE}. We state and
prove a very general error bound and compare it with classical results
by Benssousan and Lions. In Section \ref{Sec:MonS} we obtain
error bounds for equations that are non-convex in the 0-th order
term. These results are of auxiliary nature and are needed in
Section \ref{Sec:ObsProbl}. In this section we state and prove our
main result, an error bound for \eqref{S}. Then in Section
\ref{Sec:App}, we apply our main result to obtain error bounds for a
FDM and a control scheme. Finally there is an Appendix containing some
technical a priori estimates.

\section{The penalization method}
\label{Sec:PM}

In this section we will use comparison arguments to derive new error
bounds for the convergence of the solution of the penalization problem
\eqref{PE} to the solution of the obstacle problem \eqref{E}. 
Here we will no longer assume \eqref{defF}, in stead we will allow for very
general structure of $F$:\\

\noindent (C1) (Comparison) The equations \eqref{E} and \eqref{PE}
satisfy the comparison principle for the class of weak solutions
under consideration.\\ 

\noindent (C2) (Monotonicity) Let
$X,Y\in\BS^N,p,x\in\R^N,r,s\in\R$. If $X\geq Y$ and $r\leq s$, then
$$F(x,r,p,X)\leq F(x,s,p,Y).$$
\vspace*{0cm}

\noindent (C3) (Regularity) One of the following statements hold:
\begin{itemize}
\item [(i)] $g\in  C^{0,1}(\R^N)$, $|D^2g^-| \leq C$, and for every
    $x\in\R^N$ and $\phi\in C^2(\R^N)$ satisfying
$|\phi|_{0,1}+|D^2\phi^-|_{0}\leq R$, 
$$F(x,\phi(x),D\phi(x),D^2\phi(x))\leq C_R.$$
\item [(ii)] $g\in C^{0,1}(\R^N)$, and for every
    $x\in\R^N$ and $\phi\in C^2(\R^N)$ satisfying $|\phi|_{0,1}\leq R$,
$$F(x,\phi(x),D\phi(x),D^2\phi(x))\leq C_R(1+|D^2\phi^-|_0).$$  
\item [(iii)] $g\in C^{0,\mu}(\R^N)$ for some $\mu\in(0,1)$, and for every
    $x\in\R^N$ and $\phi\in C^2(\R^N)$ satisfying $|\phi|_{0}\leq R$,
$$F(x,\phi(x),D\phi(x),D^2\phi(x))\leq C_R(1+|D\phi|_0+|D^2\phi^-|_0).$$
\end{itemize}

Assumption (C1) is not very precise. In applications we need to
specify both the notion of weak solutions and ``boundary conditions''
at infinity. 
Assumption (C2) says that $F$ is ``proper'' in the terminology of the
User's Guide \cite{CIL:UG}, and implies that $F$ is degenerate elliptic.
Assumption (C3) gives regularity assumption on the obstacle
$g$ and corresponding (local) boundedness assumptions on $F$. 
Assumption (C3) can be generalized to
allow for super-linear growth in $|X^-|$ and $|p|$. This would affect
the rates obtained and will not be considered here.

These assumptions are satisfied by a {\em very}
wide class equations and with different concepts of weak solutions. In the
viscosity solutions setting (the weakest notion allowed here), we will just
mention that the above assumption hold for the Bellman equations from
stochastic control \cite{FS:Book} and the Isaacs equations from stochastic
differential games \cite{FS:SDG} under natural assumptions on the
data. We refer to the User's Guide \cite{CIL:UG} for many more
viscosity solution examples. Typical ``boundary conditions'' would be
to assume bounded solutions or linear growth at infinity. 
We can also consider variational solutions \cite{BL:Book} whenever it
makes sense to do so. In this case all point-wise inequalities have to be
interpreted in the almost everywhere sense.

The main result in this section gives both the
convergence and the rate of convergence for the penalization problem.

\begin{theorem}
\label{main}
Assume (C1) -- (C3) hold and $u$ and $v_\eps$ are solutions
of \eqref{E} and \eqref{PE} (we do not assume \eqref{defF}!). Then if
$|D^2g^-|_0<\infty$ (case (i))
$$0\leq u-v_\eps\leq C\eps \quad \text{in}\quad \R^N,$$
otherwise ($g\in C^{0,\mu}(\R^N)$ -- cases (ii) and (iii))
$$0\leq u-v_\eps\leq C\eps^{\mu/2} \quad \text{in}\quad \R^N,$$
where the constants $C$ only depend on $g$ and $C_R$ from (C3).
\end{theorem}

\begin{remark}
The rates depends only on the regularity of the obstacle, and not 
on the regularity of the solution. Even if the solution $u$ is only
H{\"o}lder continuous, we still get rate $1$ if $|D^2g^-|_0<\infty$.
For many other types of approximation schemes the rates depends
directly on the regularity of the solution, see e.g. \cite{BJ:Rate} (FDMs)
and \cite{JK:Ell} (vanishing viscosity method).  
\end{remark}

To the best of the our knowledge, this is the first time the
penalization error has been estimated for degenerate equations, and
the above result seem to be new even in the linear uniformly elliptic
case (see below). 

Before giving the proof of Theorem \ref{main}, let us briefly consider the
linear uniformly elliptic case. Here $L^2/W^{1,2}$-estimates on the
penalization error are classical \cite{BL:Book}. We will  
state a typical such result, so that the reader can compare it with
the one we have obtained. In our notation:   
\begin{align}
\label{BLOP}
\min\left\{-A u +f(x); u-g(x)\right\}&=0 && \text{in}\quad
\Omega,\\
\label{BLPP}
-A v_\eps +f(x)&=\frac{1}{\eps}(v_\eps-g(x))^- &&\text{in}\quad \Omega,%\\
\end{align}
where $\Omega$ is a smooth bounded domain and $A$ is a
linear elliptic operator in divergence form,
\begin{align*}
A \phi(x):= \del_{x_i}( a_{ij}(x)\del_{x_j}\phi)+b_i(x)\del_i\phi-\lambda \phi.
\end{align*}
The summation convention is used, $\lambda>0$, and ellipticity means 
$\xi_ia_{ij}(x)\xi_j\geq \alp|\xi|^2$ for some 
$\alp>0$ and every $\xi\in\R^N$. 
The concept of solutions is that of variational (weak)
solutions belonging to $W^{1,2}_0(\Omega)$, see \cite{BL:Book} for the
exact definitions. Typical assumptions on the data are\\

\noindent (D) $\quad a_{ij}\in L^{\infty}(\Omega), b_i\in W^{1,\infty}(\Omega),
f\in L^2(\Omega), 
g\in W^{1,2}(\Omega),\ 
Ag\in L^2(\Omega),$\\

\noindent where $a=(a_{ij})_{ij}, b=(b_i)_i$, and the error bound
obtained is the following \cite[p. 197]{BL:Book}:
\begin{proposition}[Bensoussan \& Lions] 
\label{BL_thm}
Assume (D) holds, $\lambda>0$ large enough, $g|_{\del\Omega}\geq 0$, and 
 $u$  and $v_\eps$ solve \eqref{BLOP} and \eqref{BLPP}. Then
$$\|u-v_\eps\|_{W^{1,2}(\Omega)}\leq C \eps^{1/2}.$$
\end{proposition}

Note that in this theorem we need control over the second
derivatives of the obstacle $g$ ($Ag\in L^2$ essentially means that 
$g\in W^{2,2}$), while in our result we only need to control the first
derivative of $g$ (say $g\in W^{1,\infty}=C^{0,1}$). Furthermore, we
may use Theorem \ref{main} to get a new error bound in
this case. Comparison principles for \eqref{BLOP} and \eqref{BLPP} are 
essentially given by Theorems 1.2 and 1.4 p. 192 and p. 198 in
\cite{BL:Book},  
so if $g\in W^{1,\infty}(\Omega)$ we can
conclude by Theorem \ref{main} that 
$$\|u-v_\eps\|_{L^\infty(\Omega)}\leq C \eps^{1/2}.$$
Here $u$ and $v_\eps$ are variational solutions of \eqref{BLOP} and
\eqref{BLPP}. This result does not follow from Proposition \ref{BL_thm}
unless $\Omega$ is a domain in $\R^1$.

\subsection*{The proof of Theorem \ref{main}}

We give a series of simple lemmas that leads the way to the proof of
Theorem \ref{main}. We start by a preliminary error estimate:

\begin{lemma}
\label{uv-err}
Assume (C1) and (C2). Let $u$ and $v_\eps$ solve \eqref{E} and
\eqref{PE}. Then 
$$0\leq u-v_\eps\leq |(v_\eps-g)^-|_0\quad\text{in}\quad\R^N.$$
\end{lemma}

\begin{proof}
First we check that by monotonicity in $r$ (C2),
$$v_\eps+|(v_\eps-g)^-|_0$$ 
is a supersolution of
\eqref{E}. The comparison principle for \eqref{E} then yields
the second inequality. Similarly, the first inequality follows since
$v_\eps$ is subsolution of \eqref{E}.
\end{proof}

Now we will estimate $|(v_\eps-g)^-|_0$: 

\begin{lemma}
\label{vg-err}
Assume (C1) -- (C3) hold and $g\in C^2(\R^N)$. 
According to (C3) define:
\begin{itemize}
\item [] Case (i): $\ \ K:=C_R$ with $R=|g|_0+|Dg|_0+|D^2g^-|_0$.
\item [] Case (ii): $\ K:=C_R(1+|D^2g^-|_0)$ with $R=|g|_0+|Dg|_0$. 
\item [] Case (iii): $K:=C_R(1+|Dg|_0+|D^2g^-|_0)$ with $R=|g|_0$. 
\end{itemize}
Let $v_\eps$ be the solution of \eqref{PE}. Then
\begin{align*}
&-\eps K\leq v_\eps-g\quad\text{in}\quad\R^N.
\end{align*}
\end{lemma}
\begin{proof}
The result follows from the comparison principle since 
$$g-\eps K$$ 
is a (classical) subsolution of \eqref{PE}.
\end{proof}

Since we did not assume that $g$ is smooth, we need an
approximation result. Let $\rho_\delta$ be the standard 
mollifier, $\rho_\delta(x)=1/\delta^N\rho(x/\delta)$ where $\rho$
is a smooth positive function with mass one and support in the unit
ball. Let $g_\delta=g*\rho_\delta$ and denote by $u^\delta$ and
$v^\delta_\eps$ the solutions of \eqref{E} and \eqref{PE} when
$g_\delta$ has replaced $g$:  
\begin{align}
\label{SmE}
\min\{F(x,u^\delta,Du^\delta,D^2u^\delta);u^\delta-g_\delta\}&=0
&& \text{in}\quad \R^N,\\
\label{SmPE}
F(x,v^\delta_\eps,Dv^\delta_\eps,D^2v^\delta_\eps)&=\frac{1}{\eps}(v^\delta_\eps-g_\delta)^- && \text{in}\quad \R^N.
\end{align}

We have the following bounds on $u-u^\delta$ and $v_\eps-v_\eps^\delta$:
\begin{lemma}
\label{approx1}
Assume (C1) and (C2), and let $u$, $u^\delta$, $v_\eps$ and
$v^\delta_\eps$ be solutions of \eqref{E}, \eqref{SmE},
\eqref{PE}, and \eqref{SmPE}. Then
$$|u-u^\delta|_0+|v_\eps-v^\delta_\eps|_0\leq |g-g^\delta|_0.$$
\end{lemma}

\begin{proof}
We only prove the $v$-result, the proof of the $u$-result is
similar. (If we knew a priori that $v_\eps\ra u$, the 
$u$-result could be obtained by going to the limit in the
$v$-result.) Let $K:=|g-g^\delta|_0$ and define  
$$w^\pm(x)=v_\eps^\delta(x)\pm K.$$
The result follows by the comparison principle for \eqref{PE} since
$w^+$ and $w^-$ are super- and subsolutions of \eqref{PE} respectively.

Let us prove that $w^+$ is a supersolution of \eqref{PE}, the
subsolution part is similar. First observe that by (C2)
\begin{align*}
F(x,w^+,Dw^+,D^2w^+) &\geq
F(x,v_\eps^\delta,Dv_\eps^\delta,D^2v_\eps^\delta).
\end{align*}
Then observe that by the definition of $K$,
\begin{align*}
-(w^+-g)^-=-(v_\eps^\delta+ K-g)^-\geq -(v_\eps^\delta-g_\delta)^-.
\end{align*} 
Since $v_\eps^\delta$ is a supersolution
of \eqref{SmPE}, the above observations show (at least formally) that
\begin{align*}
&F(x,w^+,Dw^+,D^2w^+)-\frac{1}{\eps}(w^+-g)^-\\
&\geq
F(x,v_\eps^\delta,Dv_\eps^\delta,D^2v_\eps^\delta)-\frac{1}{\eps}(v_\eps^\delta-g_\delta)^-\geq
0. 
\end{align*}
The proof is complete since all the above computations easily can be
seen to hold in the weak/viscosity sense.
\end{proof}

Now we can give the proof of Theorem \ref{main}:

\begin{proof}[Proof of Theorem \ref{main}]
First we consider the solutions $u^\delta$ and $v_\eps^\delta$ of
\eqref{SmE} and \eqref{SmPE}. By Lemmas \ref{uv-err} and
\ref{vg-err} we have 
$$u^\delta-v^\delta_\eps\leq K\eps \quad \text{in}\quad \R^N,$$ 
where $K$ is defined in  Lemma \ref{uv-err}. Since 
$$u-v_\eps=
(u-u^\delta)+(u^\delta-v_\eps^\delta)+(v_\eps^\delta-v_\eps),$$
Lemma \ref{approx1} and H{\"o}lder continuity of $g$ lead to
$$u-v_\eps\leq 2[g]_\mu\delta^\mu + K\eps.$$
If $|g_{xx}^-|_0<\infty$ then $K$ is independent of $\delta$ and we
can send $\delta\ra0$, leading to
$$u-v_\eps\leq K\eps.$$
Otherwise, by H{\"o}lder continuity of $g$, $K=C \delta^{\mu-2}$, and
minimization w.r.t. $\delta$ yields
$$u-v_\eps\leq C \eps^{\mu/2}.$$
The lower bound follow from Lemma \ref{uv-err}. 
\end{proof}

\begin{remark}
\label{SSec:Rem}
The procedure used in the above proof is very general, and works for
any problem satisfying the assumptions  corresponding to (C1) 
-- (C3). E.g. one could consider boundary value problems where one
would find that the estimates in Theorem \ref{main} still holds. In the
next section we will even see the method applied to an obstacle
problem for an approximation scheme (Lemma \ref{eps-scheme}). 
\end{remark}

We end this section by indicating an alternative approach. Remember
that the lower bounds in Theorem \ref{main} follow from the comparison
principle since the solution $v_\eps$ of \eqref{PE} is a subsolution
of \eqref{E}. To obtain the upper bounds, we need in some sense to
show that $v_\eps$ is an approximate supersolution of
\eqref{E}. Observe that formally 
$$\min\left\{F[v_\eps] ; v_\eps-g\right\}\geq -\eps (F[v_\eps])^+.$$
This follows since
$$0=F[v_\eps]-\frac{1}{\eps}(v_\eps-g)^-
=\min\left\{F[v_\eps] ; F[v_\eps]+\frac{1}{\eps}(v_\eps-g)\right\}$$ 
implies that
\begin{align*}
0&=\min\left\{F[v_\eps] ; \eps F[v_\eps]+v_\eps-g\right\}
\leq  \min\left\{F[v_\eps] ; v_\eps-g\right\} + \eps (F[v_\eps])^+.
\end{align*}
This is vanishing viscosity(!), and we should already guess that the
error should be $C \eps^{1/2}$ when solutions are Lipschitz
continuous. 
An easy way to get this result is the continuous dependence approach
of \cite{JK:ContDep,JK:Ell} which leads
$$u-v_\eps\leq C{``}\eps (F[v_\eps])^+{``} = C(1+|Dv_\eps|_0)
\eps^{1/2},$$
where $u$ is the solution of \eqref{E}. 
The loss of rate is caused by $v_\eps$ being only Lipschitz
continuous while $F$ is a second order operator. On the other hand, if
$|D^2v_\eps^-|_0<\infty$ then we would have the full rate:
$$u-v_\eps\leq C{``}\eps (F[v_\eps])^+{``} =
C(1+|Dv_\eps|_0+|D^2v_\eps^-|_0) \eps.$$
In Theorem \ref{main} this last estimate is proved under the much
weaker assumption $|D^2g^-|_0<\infty$.

\section{Monotone Approximation Schemes - Preliminaries.}
\label{Sec:MonS}

In this section we will give a slight generalization of the results of
\cite{Kr:HJB1,Kr:HJB2,BJ:Rate,J:Par}. We prove error bounds for
monotone approximation schemes for 
equations that have possibly non-convex dependence on $0$-order
terms. These results will then be used in Section \ref{Sec:ObsProbl} to obtain
rates for the more difficult obstacle problem \eqref{E} and \eqref{defF}.

Consider the following equation:
\begin{align}
\label{SLin}
F(x,u,Du,D^2u)= 0 \quad\text{in}\quad\R^N,
\end{align}
where $F$ is given by \eqref{defF} in the introduction.
We make the following assumptions:\\

\noindent {\bf (A1)} $a^\alp=\frac12\sigma^{\alp}{\sigma^{\alp}}^T$ for some
$N\times P$ matrix $\sigma$, and there is a $C$ independent of $\alp$
such that $|\sigma^\alp|_1+ |b^\alp|_1 \leq C.$\\ 

\noindent {\bf (A2)}  There are $\lambda,\Lambda>0$ such that for 
every $x\in\R^N$, $\alp\in\A$, $r,s\in\R$ satisfying $r\geq s$,
$$\lambda(r-s)\leq f^\alp(x,r)-f^\alp(x,s)\leq \Lambda (r-s).$$
Furthermore, $f^\alp(\cdot,0)$ is bounded uniformly in $\alp$,
and for every $x,y\in\R^N$, $r\in\R$, and $\alp\in\A$
$$|f^\alp(x,r)-f^\alp(y,r)|\leq C(1+|r|)|x-y|.$$

\begin{remark}
The first part of assumption (A2) implies that $f$ is Lipschitz and
strictly increasing in $r$. The second part implies that $f$ is
bounded and Lipschitz in $x$ for fixed $r$. 
If (A1) and (A2) hold and $\eps>0$ is fixed, the   
penalization scheme \eqref{PE} can be rewritten in the form
\eqref{SLin} by redefining $f^\alp(x,r)$ to be
$$f^\alp(x,r) +\frac1\eps\min\left\{r-g(x);0\right\}.$$
This new function than satisfies (A2) with new constants $\lambda$,
$\Lambda+\frac1\eps$, and $C$.
\end{remark}

Existence, uniqueness, and regularity follow from standard viscosity
solutions arguments. The results parallels the one mentioned in
Section \ref{Sec:intro}, and we state them without proofs:
\begin{lemma}
\label{exuniq}
Assume (A1) and (A2) hold. Then there is a unique bounded H{\"o}lder
continuous viscosity solution $u$ of \eqref{SLin}. Furthermore, if
$\lambda$ is big enough (compared to $[\sigma]_1$ and $[b]_1$), then $u$ is
Lipschitz continuous.
\end{lemma}

Using notation from the introduction, we may write an approximation scheme
for \eqref{SLin} in the following way
\begin{align}
\label{scheme3}
S(h,x,u_h(x),[u_h]_{x})=0\quad \text{in}\quad\R^N.
\end{align}
We require $S$ to satisfy:\\ 

\noindent {\bf (S1)} (Monotonicity)
For every $h>0$, $x\in\R^N$, $r\in\R$, $m \geq 0$ and bounded
functions $u,v$ such that $u\leq v$ in $\R^N$, the following holds:
$$S(h,x,r+m,[u+m]_{x}) \geq \lambda m + S(h,x,r,[v]_{x}),
$$
where $\lambda>0$ is given by (A2).\\

\noindent {\bf (S2)} (Regularity)
For every $h>0$ and $\phi\in C_b(\R^N)$,
$x \mapsto S(h,x,\phi(x),[\phi]_{x})$
is bounded and
continuous in $\R^N$ and the function $r \mapsto
S(h,x,r,[\phi]_{x})$ is uniformly continuous for bounded $r$,
uniformly in $x \in \R^N$.\\

\noindent {\bf (S3)} (Consistency) There exists integers
$n,k_i>0$, constants $K_i\geq0$, $i=1,2,\dots,n$ such 
that for every smooth $\phi$, $h > 0$, and $x\in\R^N$:
\begin{gather*}
\left|F(x,\phi(x),D\phi(x),D^2\phi(x))-S(h,x,\phi(x),[\phi]_{x})\right|\leq
\sum_{i=1}^nK_i|D^i\phi(x)|{h}^{k_i}.  
\end{gather*}

Condition (S1) and (S2) imply a comparison result for bounded
continuous solutions of \eqref{scheme3} (cf. \cite{BJ:Rate}): 

\begin{lemma}
\label{CompLem}
Assume (S1), (S2), and $u, v\in C_b(\R^N)$. If
$S[u]\leq 0$ and $S[v]\geq 0$ in  $\R^N$,  then $u\leq v$ in $\R^N.$
\end{lemma}

We proceed with obtaining an upper bound on the error for the
scheme \eqref{scheme3}. In order to do so we will consider the
following auxiliary problem:
\begin{align}
\label{SLinAux}
\sup_{|e|\leq\delta}\tilde{F}(x+e,u^\delta(x),Du^\delta(x),D^2u^\delta(x))=0
\quad\text{in}\quad\R^N,
\end{align}
where $\delta>0$, and with $u$ being the solution of \eqref{SLin},
\begin{align*}
\tilde{F}(x,r,p,X):=\sup_{\alp\in\A}\left\{-\tr[a^{\alp}(x)X]+b^\alp(x)p
  +\lambda r-\lambda u(x) + f^\alp(x,u(x))\right\}.
\end{align*}
Actually this is a problem of the same type as \eqref{SLin} so
well-posedness follows in the same way. At this point we assume the
following:\\ 

\noindent {\bf (A3)} Let $u$ and $u^\delta$ denote the solutions of
\eqref{SLin} and \eqref{SLinAux}. There is a constant $K>0$
independent of $\delta$ such that 
$$|u^\delta|_1+\frac{1}{\delta}|u-u^\delta|_0\leq K.$$
\begin{remark}
\label{RemA3}
Assumption (A3) follows from assumptions (A1) and (A2) if $\lambda$ is
big enough. After observing that \eqref{SLin} can be written as an
Isaacs equation, this follows from Lemmas \ref{WP} and \ref{CD} in the
Appendix. In the case that $\lambda$ is not ``big enough'' things are
a little bit more complicated, we refer to \cite{BJ:Rate} for this case.
\end{remark}

Now we are in a position to derive 
an upper bound on the error for the scheme \eqref{scheme3}. 

\begin{theorem}
\label{RateSLin}
Let (A1) -- (A3), (S1) -- (S3) hold, let $u$ be the viscosity
solution of \eqref{SLin} , and let $u_h$ be a
solution of the scheme \eqref{scheme3}. Then if $h> 0$ is
sufficiently small,  
\begin{gather*}
u-u_h \leq C{h}^{\gamma}
 \quad\text{in}\quad\R^N,
\end{gather*}
where $\gamma:=\underset{i:K_i>0}{\min}\left\{\frac{k_i}{i}\right\}$ and
$C\leq \frac{K}{\lambda}(\sum_{i=1}^nK_i+2(2\lambda+\Lambda))$. 
\end{theorem}

\begin{proof} 1) We start by showing that
 $u_\delta:=\rho_\delta*u^\delta$ is a subsolution of
\begin{align}
\label{Slin2}
\tilde{F}(x,w,Dw,D^2w)=0\quad \text{in}\quad \R^N,
\end{align}
 where $\rho_\delta$ is the mollifier defined in
 Section \ref{Sec:PM}. By (A3)
\begin{gather*}
\tilde{F}(x+
e,u^{\delta}(x),Du^{\delta}(x),D^2u^{\delta}(x)) \leq 0 \quad
\text{in}\quad \R^N
\end{gather*}
for every  $|e|\leq \delta$. Hence for every $|e|\leq\delta$,
$u^{\delta}(x-e)$ is a subsolution of \eqref{Slin2}.
Then $u_\delta$ is also a subsolution of \eqref{Slin2} since it can be
viewed as the limit of convex combinations of subsolutions
$u^{\delta}(x-e)$ of the convex equation \eqref{Slin2}, we refer to
the Appendix in \cite{BJ:Rate} for the details. 

2) $u_\delta$ is an approximate subsolution to the scheme \eqref{scheme3}.
By properties of mollifiers and (A3), $u_{\delta}$ is smooth and satisfies
$$\delta^{i-1}|D^iu_{\delta}|_{0}+(2\delta)^{-1}|u-u_{\delta}|_{0}\leq K.$$
So by (A2) and the definition of $\tilde{F}$, for every $x\in\R^N$,
\begin{align*}
F(x,u^{\delta}(x),Du^{\delta}(x),D^2u^{\delta}(x))&\leq \sup_{\alp\in\A}\left| \lambda( u_\delta-u) -
  f^\alp(x,u_\delta)+ f^\alp(x,u)\right|\\
&\leq 2K(\lambda+\Lambda)\delta.
\end{align*}
Consistency (S3) then leads to
\begin{align*}
&S(h,y,u_{\delta}(y),[u_{\delta}]_{y})\leq
K\sum_{i=1}^nK_i\delta^{1-i}{h}^{k_i}+ 2K(\lambda+\Lambda)\delta =:\ol{C}. 
\end{align*}

3) By (S1), $u_\delta -\ol{C}/\lambda$ is a subsolution to the scheme
\eqref{scheme3}. By comparison, Lemma \ref{CompLem}, we have
\begin{align*}
u_\delta-u_h \leq \ol{C}/\lambda
\quad \text{in}\quad \R^N.
\end{align*}

4) Combining the above estimates yields
\begin{align*}
u-u_h= u-u_{\delta}+u_{\delta}-u_h \leq 2K\delta+\ol{C}/\lambda \quad
\text{in}\quad \R^N. 
\end{align*}
Now we can conclude by choosing 
$$\delta=\max_{i:K_i>0}\{h^{k_i/i}\}.$$
\end{proof}

\begin{remark}
\label{concave}
If we replace  $\sup$  by $\inf$ in equations
\eqref{SLin} and \eqref{SLinAux}, a similar argument would lead to a
lower bound of the error: $-C h^\gamma\leq u-u_h$.
\end{remark}

From this remark it is clear that we have the full result for
{\em semi-linear} equations (see also \cite{J:Par} for the linear case):

\begin{corollary}
\label{RateLin}
Assume \eqref{SLin} is semi-linear, i.e. that $\A$ is a singleton.
Let (A1) -- (A3), (S1) -- (S3) hold, let $u$ be the viscosity
solution of \eqref{SLin}, and let $u_h$ be a
solution of the scheme \eqref{scheme3}. Then if $h> 0$ is
sufficiently small,  
\begin{gather*}
|u-u_h|_0 \leq C{h}^{\gamma},
\end{gather*}
where $\gamma$ and $C$ are defined in Theorem \ref{RateSLin}.
\end{corollary}

Following the ideas in \cite{BJ:Rate,Kr:HJB1}, we proceed to have
obtain the full result in for more general situations.
Let $\tilde{S}$ denote the scheme $S$ when it is applied to equation
$\tilde{F}[w]=0$ where $\tilde{F}$ is defined just after
\eqref{SLinAux}, and consider 
\begin{align}
\label{SAux}
\sup_{|e|\leq\delta}\tilde{S}(h,x+e,u_h^\delta(x),[u_h^\delta]_x)=0\quad\text{in}\quad \R^N,
\end{align}
and the assumption analogous to (A3):\\

\noindent {\bf (S4)} Assume $u_h$ and $u^\delta_h$ are 
solutions of \eqref{scheme3} and \eqref{SAux}, and there is a constant $K'>0$
independent of $\delta$ such that 
$$|u^\delta_h|_1+\frac{1}{\delta}|u_h-u_h^\delta|_0\leq K'.$$
In addition we need the following assumptions of $S$:\\

\noindent {\bf(S5)} (Convexity) For any $v\in C^{0,1}(\R^N)$, $h>0$,
and $x\in\R^N$ 
$$\int_{\R^N}S(h,x,v(x-e),[v(\cdot-e)]_x)\rho
_{\delta}(e)de \geq S(h,x,(v*\rho_{\delta})(x), [v*\rho_{\delta}]_x).$$

\noindent {\bf(S6)} (Commutation with translations) For any $h>0$ small enough,
$0\leq \delta\leq 1$, $y\in \R^N$, $t\in \R$, $v\in C_{b}(\R^N)$
and $|e|\leq \delta$, we have 
$$S(h,y,t,[v]^h_{y-e}) = S(h,y,t,[v(\cdot -e)]^h_{y}).$$

\begin{remark}
While (S5) and (S6) are not very restrictive, (S4) is. 
This assumption is satisfied for control schemes in general
\cite{BJ:Rate} and for FDMs when the coefficients multiplying second
order derivatives are constants (Section  \ref{Sec:App}).
Note that (S4) is not assumed in Corollary \ref{RateLin}.
\end{remark}

It is clear that by repeating the arguments in the proof of Theorem
\ref{RateSLin}, with the schemes \eqref{scheme3} and \eqref{SAux}
taking the role of the equations \eqref{SLin} and \eqref{SLinAux}, we
obtain a lower bound on the error $-Ch^\gamma\leq u_h-u$. We
refer to \cite{BJ:Rate} for more details. Combining this result 
with Theorem \ref{RateSLin} then yields the main result in this section.

\begin{theorem}
\label{FullRateSLin}
Let (A1) -- (A3), (S1) -- (S6) hold, let $u$ be the viscosity
solution of \eqref{SLin} , and let $u_h$ be a
solution of the scheme \eqref{scheme3}. Then if $h> 0$ is
sufficiently small,  
\begin{gather*}
|u-u_h|_0 \leq C{h}^{\gamma},
\end{gather*}
where $\gamma$ is defined in Theorem \ref{RateSLin} and
$C\leq \frac{K\vee
  K'}{\lambda}(\sum_{i=1}^nK_i+2(2\lambda+\Lambda))$.  
\end{theorem}

The results in this section generalize slightly the results in
\cite{Kr:HJB1,Kr:HJB2,BJ:Rate,J:Par} which consider pure convex or
concave equations. Here we allow non-convexity (non-concavity) in the
0-th order terms.

\section{Monotone Approximation Schemes - The Main Result.}
\label{Sec:ObsProbl}

In this section we will see how to use the results of the previous
two sections to obtain error bounds for monotone
schemes \eqref{S} for the non-convex problem \eqref{E}. In Section
\ref{Sec:App} we give examples of such schemes. We assume
that $S$ in \eqref{S} satisfies assumptions (S1) -- (S3) of
Section \ref{Sec:MonS}.

First we consider the penalization problem corresponding to \eqref{S},
namely problem \eqref{PS} in the introduction.
This scheme is also an approximation scheme for the
penalization problem \eqref{PE}. Note that
\eqref{S} and \eqref{PS} themselves satisfy assumptions (S1) and
(S2) (when $S$ is appropriately redefined) and hence the comparison
principle, Lemma \ref{CompLem}, holds also for these schemes. 

We start by obtaining the rate of convergence for $v_h^\eps\ra
u_h$. It is not difficult to see that this result is a consequence of
the procedure given in Section \ref{Sec:PM} if we can prove that
assumptions corresponding to (C1) -- (C3) hold for $S$. Because
(S1) and (S2) imply comparison, they already imply assumptions
corresponding to (C1) and (C2). But it turns out that (S3) 
is not sufficient to have the assumption corresponding to (C3) because
it involves derivatives of higher order than two. We need to assume
``(C3)'':\\  

\noindent {\bf(S7)} (Regularity) One of the following statements hold:
\begin{itemize}
\item [(i)] $g\in C^{0,1}(\R^N)$, $|D^2g^-|_0 \leq C$, and for every
    $x\in\R^N$ and $\phi\in C^2(\R^N)$ satisfying
$|\phi|_{0,1}+|D^2\phi^-|_0\leq R$, 
$$S(h,x,\phi(x),\phi)\leq C_R.$$
\item [(ii)] $g\in C^{0,1}(\R^N)$, and for every
    $x\in\R^N$ and $\phi\in C^2(\R^N)$ satisfying $|\phi|_{0,1}\leq R$,
$$S(h,x,\phi(x),\phi)\leq C_R(1+|D^2\phi^-|_0).$$  
\item [(iii)] $g\in C^{0,\mu}(\R^N)$ for some $\mu\in(0,1)$, and for every
    $x\in\R^N$ and $\phi\in C^2(\R^N)$ satisfying $|\phi|_{0}\leq R$,
$$S(h,x,\phi(x),\phi)\leq C_R(1+|D\phi|_0+|D^2\phi^-|_0).$$
\end{itemize}
\begin{remark}
This assumption holds for most reasonable schemes \eqref{S} when
(S3) also holds, e.g. for the finite 
difference method \eqref{FDM} below. 
\end{remark}

By the method of Section \ref{Sec:PM}, we have the following result:
\begin{lemma}
\label{eps-scheme}
Assume (S1), (S2), (S7) hold and $u_h$ and $v_{h,\eps}$ are solutions
of \eqref{S} and \eqref{PS}. Then if $|D^2g^-|_0<\infty$
(case (i))
$$0\leq u_h-v_{h,\eps}\leq C\eps \quad \text{in}\quad \R^N,$$
Otherwise ($g\in C^{0,\mu}(\R^N)$ -- cases (ii) and (iii))
$$0\leq u_h-v_{h,\eps}\leq C\eps^{\mu/2} \quad \text{in}\quad \R^N,$$
where the constants $C$ only depend on $g$ and $C_R$ from (C3).
\end{lemma}

Now to obtain results for the scheme \eqref{S}, we may use the
following diagram:
\begin{equation*}
    \xymatrix{\min\{F[u]; u-g\}=0
        \ar@{<->}[rrrr]^{?}\ar@{<->}[dd]_{0\leq u-v_\eps\leq C_1(\eps)}^{\text{Theorem
        \ref{main}}}  
        & 
        &
        &
        & \min\{ S[u_h]; u-g\}=0 \ar@{<->}[dd]^{0\leq u_h-v_{h,\eps}\leq
        C_4(\eps)}_{\text{Lemma \ref{eps-scheme}}} \\
        &&&&\\
        F[v_\eps]=\frac{1}{\eps}(v_\eps-g)^- \ar@{<->}[rrrr]_{C_2(h,\eps)\leq v_\eps-v_{h,\eps}\leq C_3(h,\eps)}^{\text{Theorem \ref{FullRateSLin}}} 
        &
        &
        &  
        & S[v_{h,\eps}]=\frac{1}{\eps}(v_{h,\eps}-g)^- }
\end{equation*}
The main result of this paper is the following:

\begin{theorem}
\label{FullRateOSLin}
Let (A1), (A2), (S1) -- (S7) hold with
  $\lambda>$ $\sup_{\alp}([\sigma^{\alp}]_1^2+[b^{\alp}]_1)$
  in (A2) and $K'$ independent of $\eps$ in (S4), let $u$ be the 
viscosity solution of \eqref{E} with $F$ defined in \eqref{defF}, and
let $u_h$ be a solution of the scheme \eqref{S}. Then if $h> 0$ is
sufficiently small,  
\begin{gather*}
|u-u_h|_0 \leq C{h}^{\gamma/3}. 
\end{gather*}
If in addition $|D^2g^-|_0<\infty$, then
$$|u-u_h|_0 \leq C{h}^{\gamma/2}. 
$$
Here $\gamma$ is defined in Theorem \ref{RateSLin} and the constants
$C$ are independent of $h$.
\end{theorem}

The assumption on $\lambda$ may be relaxed to simply requiring
$\lambda>0$. This will influence the rates and complicate the
arguments, see \cite{BJ:Rate} for a discussion. See also Remark \ref{RemA3}.

\begin{proof}[Outline of proof]
1) By Lemmas \ref{WP} and \ref{CD} in the Appendix (see Remark
   \ref{RemA3}), assumption (A3) is satisfied for \eqref{PE} with
   $K$ independent of $\eps$! Note that $K'$ is assumed independent of $\eps$.

2) By Theorem \ref{FullRateSLin} with $\Lambda$ replaced by
   $\Lambda+\frac1\eps$, 
$$|v_\eps-v_{h,\eps}|_0\leq C(1+\frac1\eps)h^\gamma \quad \text{in}\quad\R^N,$$
where $v_{h,\eps}$ solves \eqref{PS} and $C$ is a constant independent
of $\eps$.

3) The result now follows from the triangle
   inequality, part 2), Theorem \ref{main}, Lemma \ref{eps-scheme}, and a
   minimization in $\eps$ (see the above diagram).
\end{proof}

If $F$ is concave instead of convex so that the obstacle problem
\eqref{E} is concave, then we obtain better rates using directly
Theorem \ref{FullRateSLin}: 
$$|u-u_h|_0\leq Ch^\gamma.$$
This was essentially the case considered by
\cite{BJ:Rate,J:Par}. Theorem \ref{FullRateOSLin} is the first result
for multi-dimensional non-concave/non-convex equation.

\section{Applications}
\label{Sec:App}

\subsection{A finite difference scheme}

In this section we apply a finite difference scheme proposed by
Kushner \cite{KD:Book} to the $N$-dimensional non-convex equation \eqref{E}
where $F$ is given by \eqref{SLin} and the coefficient $a$ is
independent of $x$.

We will assume that (A1) and (A2) of Section \ref{Sec:MonS} and
that the following assumptions hold:\\  

\noindent  {\bf (A4)} $\quad a^{\alp}$ is independent of $x$,\\

\noindent  {\bf (A5)} $\quad a_{ii}^{\alp} - \sum_{j\neq i}
|a_{ij}^{\alp}| \geq 0, \quad i=1,\dots,N,$\\

\noindent  {\bf (A6)} $\quad
\sum_{i=1}^{N}\Big\{a_{ii}^{\alp} - 
\sum_{j\neq i} |a^{\alp}_{ij}| + |b^{\alp}_i(x)|
\Big\}\leq 1 \quad \hbox{   in }\R^N.$\\

\noindent  {\bf (A7)} $\quad$ $\sup_{\alp}\left\{\inf_x
c^{\alp}-2\sqrt{N}[b^{\alp}]_1\right\}=:\lambda_0>0.$\\

\noindent  {\bf (A8)} $\quad$ (i) $g\in C^{0,1}(\R^N)$ or (ii) $g\in
C^{0,1}(\R^N)$ and $|D^2g^-|_0 \leq C$.\\

Here we need (A4) in order to prove condition (S4) of Section
\ref{Sec:MonS}, for more on this see \cite{BJ:Rate}. To avoid (A4) we
must use the much more difficult methods of 
\cite{BJ:Rate2} or \cite{Kr:HJB2}. We will not consider this here. 
Condition (A5) simply says that $a$ is diagonally dominant. This
is a standard condition \cite{KD:Book} and implies that the
scheme \eqref{S1} below is monotone. Conditions (A6) is a
normalization of the coefficients in \eqref{I}.  We can always have
this assumption satisfied by multiplying equation \eqref{I} by an
appropriate positive constant. Conditions (A7) and (A8) together with
(A1) and (A2) assure that the solutions of the various schemes
(e.g. \eqref{S}) belong to $C^{0,1}(\R^N)$. Under these assumptions
the solutions of various equations (e.g. \eqref{E}) will also belong
to $C^{0,1}(\R^N)$. We refer to the Appendix for the proof of these
facts. Condition (A8) is a regularity condition on $g$,
cf. (C3) and (S7).

The difference operators we use are defined in the following
way
\begin{align*}
&\Delta^{\pm}_{x_i} w(x)=\pm\frac{1}{h}\{w(x \pm e_ih)-w(x) \},\\
&\Delta^2_{x_i} w(x)= \frac{1}{h^2}\{w(x+e_ih) -2w(x) + w(x-e_ih)\},\\
&\Delta^+_{x_ix_j} w(x)= \frac{1}{2h^2}\{2w(x) + w(x+e_ih+e_jh)
+ w(x-e_ih-e_jh)\} \\ &\qquad \qquad \quad
- \frac{1}{2h^2}\{w(x+e_ih) + w(x-e_ih) + w(x+e_jh) +
w(x-e_jh)\},\\
&\Delta^-_{x_ix_j} w(x)= \frac{1}{2h^2}\{w(x+e_ih) + w(x-e_ih) +
w(x+e_jh) + w(x-e_jh)\} \\ &\qquad \qquad \quad-
\frac{1}{2h^2}\{2w(x) + w(x+e_ih-e_jh) + w(x-e_ih+e_jh)\}.
\end{align*}
Let $b^+=\max\{b,0\}$ and $b^-=(-b)^+$. Note that $b=b^+-b^-$. For each
$x$, $t$, $p^{\pm}_i$, $A_{ii}$, $A_{ij}^{\pm}$, $i,j=1,\dots,N$, let
\begin{align*}
&\tilde{F}(x,r,p^{\pm}_i,A_{ii},A_{ij}^{\pm})\\
&=\min\bigg\{\sup_{\alp\in\A} \Big\{\sum_{i=1}^N
\Big[-\frac{a^{\alp}_{ii}}{2} A_{ii} + \sum_{j\neq i}
\Big( -\frac{a^{\alp+}_{ij}}{2} A_{ij}^++
\frac{a^{\alp-}_{ij}}{2} A_{ij}^- \Big)\\
&\qquad- b_i^{\alp +}(x) p^+_i + b_i^{\alp -}(x)
p^-_i \Big] + f^{\alp}(x,r)
\Big\},r-g(x)\bigg\}.
\end{align*}
Now we can write the finite difference scheme in the following way,
\begin{align}
\label{FDM}
\tilde{F}(x,u_h(x),\Delta^{\pm}_{x_i}u_h(x),\Delta^2_{x_i}u_h(x),\Delta^{\pm}_{x_ix_j}u_h(x))=0.
\end{align}
This is a consistent and monotone scheme. 

In order to get our result, we must define $S$ in \eqref{S} and prove
that conditions (A3), (S1) -- (S7) of Sections \ref{Sec:MonS} and
\ref{Sec:ObsProbl} hold. We have moved most of the details to Appendix
\ref{App:FDM} where a more general problem is considered. 
To see how $S$ may be defined, see \eqref{defS} in Appendix
\ref{App:FDM}. Condition (S1) holds
by monotonicity of the scheme, (S2) holds trivially, (S3) holds with
following estimate: 
$$
|F(x,v,Dv,D^2v)-S(h,x,v(x),[v]^h_x)| \leq \bar{K} (|D^2v|_0 h+|D^4v|_0h^2),
$$
for any $v\in C^4(\R^N)$. Condition (S5) holds by ``convexity'' of the
$\sup$-part of the scheme, (S6) holds trivially, and by (A8) we
immediately get (S7). The only difficult
condition is (S4). To prove it we need very precise a priori estimates
on the scheme 
provided by Lemmas \ref{u_h-ex} and \ref{u-u} in Appendix
\ref{App:FDM}. Note in particular that the bounds in (S4) are
independent of the penalization parameter $\eps$.

In view of Theorem \ref{FullRateSLin} we have the following result:
\begin{proposition}
\label{LinSRate}
Assume (A1),(A2), (A4) -- (A8) of Sections \ref{Sec:MonS} and
\ref{Sec:App} hold, $u$ is the viscosity solution of \eqref{E}
and $u_h$ is the solution of \eqref{FDM}. Then if $h>0$ is
sufficiently small,
$$|u-u_h|_0\leq C h^{1/6} \quad\text{in}\quad \R^N.$$
If in addition $|D^2g^-|_0<\infty$, then
$$|u-u_h|_0\leq C h^{1/4} \quad\text{in}\quad \R^N.$$
\end{proposition}

In the convex case under similar assumptions the rate is $1/2$ when
$a$ is independent of $x$ \cite{J:Par} and at least $1/5$ in the
general case \cite{BJ:Rate2}. For one-dimensional non-convex problems
the rate is at least $1/5$ \cite{J:Non-conv}, and for first order
problems the rate is again $1/2$ \cite{CL:Appr}.

\subsection{Control schemes}
In this section, we consider a so-called control schemes introduced in
the second order case by Menaldi \cite{Men}. The scheme is defined in the
following way, 
\begin{gather}
\label{c-scheme}
u_h(x) = \min_{\vartheta\in\Theta} \Big\{ (1-h c^{\vartheta}(x))
\Pi_h^{\vartheta} u_h(x) +hf^{\vartheta}(x) \Big\}, \\
\intertext{where $\Pi_h^{\vartheta}$ is the operator defined by}
\begin{split}
&\Pi_h^{\vartheta} \phi(x) = \\
&\frac{1}{2N} \sum_{m=1}^N \Big(
\phi(x+hb^{\vartheta}(x)+\sqrt{h}\sigma^{\vartheta}_m(x))+\phi(x+hb^{\vartheta}
(x)-\sqrt{h}\sigma^{\vartheta}_m(x)) \Big),\nonumber
\end{split}
\end{gather}
and
$\sigma_m^{\vartheta}$ is the $m$-th column of $\sigma^{\vartheta}$.
In the convex case a fully discrete method is derived from
\eqref{c-scheme} and analyzed in \cite{CF:Appr}. The authors also
provide an error bound for the convergence of the solution of the
fully discrete method to the solution of the scheme \eqref{c-scheme}.

In this case we only need to assume conditions (A1), (A2), (A7), and
(A8), in particular $a$ may depend on $x$. All condition (S1) -- (S8)
then holds, and the consistency condition (S4) takes the form
$$|F(x,v,{D}v,{D}^2v)-S(h,x,v(x),[v]^h_x)| \leq \bar{K}
(|D^2v|_0+|D^3v|_0+|D^4v|_0) h,
$$
for any $v\in C^{4}(\R^N)$. We refer to \cite{BJ:Rate} for the proof
of these conditions and the precise definition of $S$. The only
difficult point is again (S4). To prove this condition one must modify the
arguments of \cite{BJ:Rate} in a similar way to what we did in the
Appendix for the FDM. We omit the details.

In view of Theorem \ref{FullRateSLin}, we have the following result:
\begin{proposition}
\label{ctrl}
Assume (A1),(A2), (A7), and (A8) of Sections \ref{Sec:MonS} and
\ref{Sec:App} hold, $u$ is the viscosity solution of \eqref{E}
and $u_h$ is the solution of \eqref{c-scheme}. Then if $h>0$ is
sufficiently small,
$$|u-u_h|_0\leq C h^{1/12} \quad\text{in}\quad \R^N.$$
If in addition $|D^2g^-|_0<\infty$, then
$$|u-u_h|_0\leq C h^{1/8} \quad\text{in}\quad \R^N.$$
\end{proposition}

In the convex case under similar assumptions the rate is at least $1/4$
\cite{J:Par}, if the solution in addition has 3 bounded derivatives then
the rate is at least $1/2$ \cite{Men}. For one-dimensional non-convex
problems the rate is at least $1/10$ \cite{J:Non-conv}, and for first order
problems the rate is again $1/2$, see e.g. \cite{BJ:Rate}.

\appendix
\section{Estimates on the Isaacs equation.}
\label{App:I}

In this section we will give well posedness results and very precise a
priori bounds for the Isaacs equation 
\begin{align}
\label{I}
\inf_{\alp\in\A}\sup_{\beta\in\B}\left\{-\tr\left[a^{\alp,\beta}(x)D^2u\right]-
b^{\alp,\beta}(x)Du+c^{\alp,\beta}(x)u-f^{\alp,\beta}(x)\right\}=0
\end{align}
in $\R^N$, where $a=\sigma \sigma^T$ for some matrix (function)
$\sigma$. We take the following assumption:\\

\noindent {\bf (B1)} $c>0$ and there is a constant $C$ independent of
$\alp,\beta$ such that 
$$[\sigma^{\alp,\beta}]_1+ [b^{\alp,\beta}]_1
+[c^{\alp,\beta}]_1+|f^{\alp,\beta}|_1\leq C.$$

We start by existence, uniqueness, and $L^\infty$-bounds on the
solution and its gradient.
\begin{lemma}
\label{WP}
If (B1) holds and
$\sup_{\alp,\beta}\left\{\inf_x c^{\alp,\beta}-[\sigma^{\alp,\beta}]_1^2-[b^{\alp,\beta}]_1\right\}>0$,   
then there exists a unique solution $u$ of \eqref{I} satisfying the
following bounds:
$$|u|_0\leq \sup_{\alp,\beta}
\frac{|f^{\alp,\beta}|_0}{\inf_x c^{\alp,\beta}},\quad
|Du|_0\leq \sup_{\alp,\beta}
\frac{|u|_0[c^{\alp,\beta}]_1+[f^{\alp,\beta}]_1}{\inf_x c^{\alp,\beta}-[\sigma^{\alp,\beta}]_1^2-[b^{\alp,\beta}]_1}.$$ 
\end{lemma}

\begin{remark}
Usually the assumption on $c$ is
$c\geq\lambda>\sup_{\alp,\beta}\left\{[\sigma^{\alp,\beta}]_1^2+[b^{\alp,\beta}]_1\right\}$
and all estimates are given in terms of $\lambda$ instead of $c$. For
our purpose this is not good enough, since we need to consider limit
problems where for some values of $\alp,\beta$, both $|f|$ and $|c|$
blow up, while for others they both remain bounded
(cf. the penalization method). 
\end{remark}

\begin{proof}
Existence and uniqueness follows from the (strong) comparison
principle and Perron's method \cite{Is:Unique}. Let 
$$M:=\sup_{\alp,\beta} 
\frac{|f^{\alp,\beta}|_0}{\inf_x c^{\alp,\beta}},$$
then the first bound on $u$ follows from the comparison principle
after checking that $M$ ($-M$) is a supersolution (subsolution) of
\eqref{I}. To get the bound on the gradient of $u$, consider
$$m:=\sup_{x,y\in\R^N}\left\{u(x)-u(y)-L|x-y|\right\}.$$
If by setting 
$$L:=\sup_{\alp,\beta}
\frac{|u|_0[c^{\alp,\beta}]_1+[f^{\alp,\beta}]_1}{\inf_x c^{\alp,\beta}-[\sigma^{\alp,\beta}]_1^2-[b^{\alp,\beta}]_1},$$
we can conclude that $m\leq 0$, then we are done. Assume for
simplicity that the maximum is attained in $(\bx,\by)$. If $\bx=\by$
then $m=0$ and we are done. If not, then $L|x-y|$ is smooth at
$(\bx,\by)$ and a standard doubling of variables argument leads to
$m\leq 0$.
Since the maximum need not be attained, we must modify the test
function in the standard way. We skip the details. (The interested
reader can have a look at the appendix of \cite{Is:Equiv} where the
above argument is given for a linear equation.)
\end{proof}

Now we proceed to obtain continuous dependence on the
coefficients. Let $\bu$ solve the following equation: 
\begin{align}
\label{bI}
\inf_{\alp\in\A}\sup_{\beta\in\B}\left\{-\tr\left[\ba^{\alp,\beta}(x)D^2\bu\right]- 
\bb^{\alp,\beta}(x)D\bu+\bc^{\alp,\beta}(x)\bu-\bff^{\alp,\beta}(x)\right\}=0
\end{align}
in $\R^N$, where $\ba=\bs \bs^T$ for some matrix (function) $\bs$. 
\begin{lemma}
\label{CD}
If $u$ and $\bu$ are bounded Lipschitz continuous solutions of
\eqref{I} and \eqref{bI} respectively, and that both sets of
coefficients satisfy (B1). Then 
\begin{align*}
&|u-\bu|_0\leq 
\sup_{\alp,\beta}\frac{K}{\inf_x c \vee \inf_x \bar c}|\sigma-\bs|_0\\ 
&
+\sup_{\alp,\beta}\frac{1}{\inf_x c \vee \inf_x \bar c}\Big\{2L|b-\bb|_0+M|c-\bc|_0+|f-\bff|_0\Big\},
\end{align*} 
where $L=[u]_1\vee[\bu]_1$, $M=|u|_0\vee|\bu|_0$, and
\begin{align*}
K^2=
32L\sup_{\alp,\beta}\Big\{&4L[\sigma]^2_1\wedge[\bs]^2_1+2L[b]_1\wedge[\bb]_1+M[c]_1\vee[\bc]_1+[f]_1\wedge[\bff]_1\Big\}. 
\end{align*} 
\end{lemma}

\begin{proof}[Outline of proof]
Define
$$m:=\sup_{x,y}\left\{u(x)-\bu(y)-\frac{1}{\delta}|x-y|^2-\eps(|x|^2+|y|^2)\right\}.$$
Then do doubling of variables using the 3 last terms in the above
expression as test-function. Using the definition of viscosity
solutions and subtracting the resulting inequalities lead to
\begin{align*}
0\leq
\sup_{\alp,\beta}\Big\{&-\tr[\ba(y)Y]+\tr[a(x)X]-\bb(y)p_x+b(x)p_y\\
&+\bc(y)\bu(y)-c(x)u(x)-\bff(y)+f(x)\Big\},
\end{align*}
where $x,y$ is the maximum point for $m$ and $(p_x,X)$, $(-p_y,Y)$ are the
elements in second order semi-jets in for $u$, $\bu$ given by the
maximum principle for semi-continuous functions \cite{CIL:UG}. Now we
note that by using Lipschitz regularity of the solutions, a standard
argument yields
$$|x-y|\leq \delta L.$$
So using Ishii's trick \cite[pp. 33,34]{Is:Unique} on the 2nd order terms,
and a few other manipulations, we get
\begin{align*} 
0\leq\sup_{\alp,\beta}\Big\{&\frac{4}{\delta}|\sigma(x)-\bs(y)|^2
+2L|b(x)-\bb(y)|+C\eps(1+|x|^2+|y|^2)\\
&+M|c(x)-\bc(y)|-(\inf_x c \vee \inf_x \bar c) m + |f(x)-\bff(y)|\Big\}.
\end{align*}
Some easy manipulations now lead to an estimate for $m$, and using the
definition of $m$, we obtain an estimate for $|u-\bu|_0$ depending on $\delta$
and $\eps$. We finish the proof by minimizing this expression
w.r.t. $\delta$ and sending $\eps\ra0$.
\end{proof}

For more details on such manipulations, we refer to \cite{JK:ContDep,JK:Ell}.

\section{Estimates on a finite difference scheme}
\label{App:FDM}

In this section we apply a finite difference scheme proposed by
Kushner \cite{KD:Book} to the $N$-dimensional Isaacs equation
\eqref{I} with coefficient $a$ independent of $x$.  

We will assume that (B1) of Section \ref{App:I} and  the following
assumptions hold:\\  

\noindent  {\bf (B2)} $\quad a$ is independent of $x$.\\

\noindent  {\bf (B3)} $\quad a_{ii}^{\alp,\beta} - \sum_{j\neq i}
|a_{ij}^{\alp,\beta}| \geq 0, \quad i=1,\dots,N,$\\

\noindent  {\bf (B4)} $\quad
\sum_{i=1}^{N}\Big\{a_{ii}^{\alp,\beta} - 
\sum_{j\neq i} |a^{\alp,\beta}_{ij}| + |b^{\alp,\beta}_i(x)|
\Big\}\leq 1 \quad \hbox{   in }\R^N.$\\

Let us define the scheme. For each
$x$, $r$, $p^{\pm}_i$, $A_{ii}$, $A_{ij}^{\pm}$, $i,j=1,\dots,N$, let
\begin{align*}
&\tilde{F}(x,r,p^{\pm}_i,A_{ii},A_{ij}^{\pm})\\
&=\inf_{\alp\in\A}\sup_{\beta\in\B} \Big\{\sum_{i=1}^N
\Big[-\frac{a^{\alp,\beta}_{ii}}{2} A_{ii} + \sum_{j\neq i}
\Big( -\frac{a^{\alp,\beta+}_{ij}}{2} A_{ij}^++
\frac{a^{\alp,\beta-}_{ij}}{2} A_{ij}^- \Big)\\
&\qquad- b_i^{\alp,\beta +}(x) p^+_i + b_i^{\alp,\beta -}(x)
p^-_i \Big] + c^{\alp,\beta}(x) r - f^{\alp,\beta}(x)
\Big\}.
\end{align*}
Using the difference operators $\Delta^{\pm}_{x_i}$, $\Delta^2_{x_i}$,
$\Delta^{\pm}_{x_ix_j}$ defined in Section \ref{Sec:App} we can now write
the finite difference scheme in the following way, 
\begin{align}
\label{S1}
\tilde{F}(x,u_h(x),\Delta^{\pm}_{x_i}u_h(x),\Delta^2_{x_i}u_h(x),\Delta^{\pm}_{x_ix_j}u_h(x))=0. 
\end{align}
This is a consistent and monotone scheme. In the following it will be
convenient to use an
equivalent formulation of this scheme (see \cite{BJ:Rate} for more
details):
\begin{align}
\label{S2}
u_h(x) = \inf_{\alp\in\A}\sup_{\beta \in\B}
\bigg\{&\frac{1}{1+h^2c^{\alp,\beta}(x)} \\
&\cdot\Big( \sum_{z\in h\Z^N} p^{\alp,\beta}(x,x+z) u_h(x+z) +
h^2f^{\alp,\beta}(x)\Big)\bigg\},\nonumber
\end{align}
where 
\begin{align*}
&p^{\alp,\beta}(x,x)=1 - \sum_{i=1}^{N}\Big\{a_{ii}^{\alp,\beta} - \sum_{j\neq i}
\frac{|a^{\alp,\beta}_{ij}|}2 + h |b^{\alp,\beta}_i(x)| \Big\},\\
&p^{\alp,\beta}(x,x\pm e_ih)= \frac{a_{ii}^{\alp,\beta}}{2} -
\sum_{j\neq i} \frac{|a^{\alp,\beta}_{ij}|}{2} + h b^{\alp,\beta\pm}_i(x),\\
&p^{\alp,\beta}(x,x+e_ih \pm e_jh)= \frac{a^{\alp,\beta \pm}_{ij}}{2},\\
&p^{\alp,\beta}(x,x-e_ih \pm e_jh)= \frac{a^{\alp,\beta \mp}_{ij}}{2},
\end{align*}
and $p^{\alp,\beta}(x,y)=0$ for all other $y$. Let $h\leq1$. Note that
by (B3) and (B4), $0\leq p^{\alp,\beta}(x,y)\leq 1$ for
all $\alp,\beta, x, y$. 
Furthermore $\sum_{z\in h\Z^N} p^{\alp,\beta}(x,x+z) = 1$ for all
$\alp,\beta, x$.

For the readers' convenience we will state explicitly the
function $S$ (as in \eqref{S} and \eqref{scheme3}) corresponding to this
scheme: We set $[\phi]^h_{x}(\cdot) := \phi(x+ \cdot)$ and 
\begin{align}
\label{defS}
&S(h,y,r,[\phi]^h_{x}):= \\ \nonumber
&\inf_{\alp\in\A}\sup_{\beta \in\B}
\left\{ - \frac{1}{h^2}\left[\sum_{z\in h\Z^N} p^{\alp,\beta}(y,y+z)
[\phi]^h_{x}(z)-t\right] +c^{\alp,\beta}(x)r - f^{\alp,\beta}(y)\right\}.
\end{align}

We use fix point arguments to prove existence, uniqueness, and a
priori bounds for equation \eqref{S1}.
\begin{lemma}
\label{u_h-ex}
Assume (B1) -- (B4) hold and 
$$\sup_{\alp,\beta}\left\{\inf_x
c^{\alp,\beta}-2\sqrt{N}[b^{\alp,\beta}]_1\right\}=:\lambda_0>0.$$
Then there exists a unique solution $u_h\in C^{0,1}(\R^N)$ of the scheme
\eqref{S1} satisfying the following bounds
$$|u_h|_0\leq \sup_{\alp,\beta}
\frac{|f^{\alp,\beta}|_0}{\inf_x c^{\alp,\beta}},\quad
[u_h]_1\leq \sup_{\alp,\beta}
\frac{\frac{|u_h|_0 +
h^2|f^{\alp,\beta}|_0}{1+h^2\inf_xc^{\alp,\beta}}[c^{\alp,\beta}]_1+[f^{\alp,\beta}]_1}{\inf_x  c^{\alp,\beta}-2\sqrt{N}[b^{\alp,\beta}]_1}.$$ 
\end{lemma}

\begin{proof}
Define $T_h:C_{b}(\R^N) \to C_{b}(\R^N)$ in the following way:
\begin{align*}
&T_{h}v(x):= \inf_{\alp\in\A}\sup_{\beta \in\B}
  \bigg\{\frac{1}{1+h^2c^{\alp,\beta}(x)}\\ 
&\qquad\qquad\qquad\qquad\cdot\Big( \sum_{z\in h\Z^N}
  p^{\alp,\beta}(x,x+z) v (x+z) + h^2f^{\alp,\beta}(x)\Big)\bigg\}.
\end{align*}
For $u,v\in C_{b}(\R^N)$, we subtract the 
expressions for $T_hu$ and $T_hv$. After we use the inequality
$\inf\sup(\cdots)-\inf\sup(\cdots)\leq \sup\sup(\cdots - \cdots)$, the
properties of $p^{\alp,\beta}$, and (B1), we obtain 
\begin{align*}
&T_hu(x)-T_hv(x) \\
&\leq  
\sup_{\alp,\beta}\,\left\{ \frac{1}{1+
  h^2\inf_xc^{\alp,\beta}}\sum_{z\in h\Z^N} 
p^{\alp,\beta}(x,x+z)|u(x+z)-v(x+z)|\right\}
\\
&\leq \frac{1}{1+\lambda_0 h^2}|u-v|_0.
\end{align*}
Since we may reverse the roles of $u$ and $v$, we see that $T_h$ is a
contraction in $(C_{b}(\R^N),|\cdot|_0)$. 
Banach's fixed point theorem then yields the existence and uniqueness
of a $u_h\in C_{b}(\R^N)$ solving \eqref{S2} (and \eqref{S1}). The
estimate on $|u_h|_0$ follows easily from the identity
$|u_h|_0=|T_hu_h|_0$. 

We proceed by proving that $u_h$ has a bounded Lipschitz
constant assuming for simplicity that $c^{\alp,\beta}$ is
independent of $x$. 
Let $v\in C^{0,1}(\R^N)$ and subtract the expressions for
$T_hv(x)$ and $T_hv (y)$:
\begin{align*}
&T_hv(x)- T_hv(y) \leq\\
&  \sup_{\alp,\beta}\bigg\{\frac{1}{1+ h^2c^{\alp,\beta}}\bigg(
\sum_{z\in h\Z^N} \Big[ p^{\alp,\beta}(x,x+z)(v(x+z)-v(y+z))\\
& + v(y+z)\big(p^{\alp,\beta}(x,x+z)-p^{\alp,\beta}(y,y+z)\big) \Big]
+ h^2 (f^{\alp,\beta}(x) -  f^{\alp,\beta}(y))\bigg) \bigg\}.
\end{align*}
In the right-hand side the first sum is bounded by $[v]_1|x-y|$, and by
using the definition of
$p^{\alp,\beta}$, the second sum is equivalent to
\begin{align*}
&h\sum_{i=1}^N\Big[\big(b_i^{\alp,\beta+}(x)-b_i^{\alp,\beta+}(y)\big)
\Delta^+_{x_i}v(y)-\big(b_i^{\alp,\beta-}(x)-b_i^{\alp,\beta-}(y)\big)\Delta
^-_{x_i}v(y)\Big]\\
& \leq  2  \sqrt{N} h^2
|b^{\alp,\beta}(x)-b^{\alp,\beta}(y)|[v]_1= 2  \sqrt{N} h^2[b^{\alp,\beta}]_1[v]_1|x-y|.
\end{align*}
Let $C^{\alp,\beta}:=2  \sqrt{N} h^2[b^{\alp,\beta}]_1$. By the above
expressions, and by exchanging the roles of $x$ and $y$, we obtain the
following estimate 
\begin{align}
\label{T-estim}
&|T_hv(x) - T_hv(y)| \leq\\
& \sup_{\alp,\beta}\bigg\{\frac{1}{1+h^2c^{\alp,\beta}}
\Big((1+h^2 C^{\alp,\beta})
    [v]_1+h^2[f^{\alp,\beta}]_1\Big)\bigg\}|x-y|.
\nonumber
\end{align}
Hence $T_hv\in C^{0,1}(\R^N)$, and $u_h\in C^{0,1}(\R^N)$ since
$u_h=\lim_{i\ra\infty}(T_h)^iv_0$ for any $v_0\in C^{0,1}(\R^N)$.
Furthermore, since $c^{\alp,\beta}\geq C^{\alp,\beta}+\lambda_0$
the  estimate on $[u_h]_1$ 
follows easily from the identity $[u_h]_1=[T_hu_h]_1$ and \eqref{T-estim}.

When $c$ depend also on $x$ we obtain an expression like
\eqref{T-estim} with $c^{\alp,\beta}$  and
$\sup_{\alp,\beta}[f^{\alp,\beta}]_1$ replaced by
$\inf_xc^{\alp,\beta}$ and
$\sup_{\alp,\beta}\big([f^{\alp,\beta}]_1 + [c^{\alp,\beta}]_1\frac{|v|_0 +
h^2|f^{\alp,\beta}|_0}{1+h^2\inf_xc^{\alp,\beta}}\big)$ respectively,
and hence the lemma would hold again.
\end{proof}

Using a standard maximum principle type of argument, we now derive a priori
estimates on the continuous dependence on the data.
\begin{lemma}
\label{u-u}
Assume (B1) -- (B4) hold and $u_h,\bar u_h\in C^{0,1}(\R^N)$. If $u_h$
solve \eqref{S1} with data $(a,b,c,f)$ and $\bar u_h$
solve \eqref{S1} with data $(a,\bar b,\bar c,\bar f)$ (same
$a$!), then 
\begin{align*}
&|u_h-\bu_h|_0\leq 
\sup_{\alp,\beta}\frac{1}{\inf_x c \vee \inf_x \bar c}\Big\{2L|b-\bb|_0+M|c-\bc|_0+|f-\bff|_0\Big\},
\end{align*} 
where $L=\sqrt{N}[u_h]_1\vee[\bu_h]_1$, $M=|u_h|_0\vee|\bu_h|_0$.
\end{lemma}

\begin{proof}
We will assume that $\sup(u-\bu)=(u-\bu)(x)\geq 0$. The general case
follows from standard modifications to the proof below. 

Since the scheme \eqref{S1} is monotone, at the maximum point $x$ we have
\begin{align*}
&\sum_{i=1}^N
\Big[-\frac{a^{\alp,\beta}_{ii}}{2} \Delta^2_{x_i} + \sum_{j\neq i}
\Big( -\frac{a^{\alp,\beta+}_{ij}}{2} \Delta_{x_ix_j}^++
\frac{a^{\alp,\beta-}_{ij}}{2} \Delta_{x_ix_j}^-
\Big)\Big](u_h-\bu_h)(x)\leq 0,
\end{align*}
and
\begin{align*}
&\sum_{i=1}^N\Big[b_i^{\alp,\beta +}(x) \Delta^+_{x_i} + b_i^{\alp,\beta -}(x)
\Delta^-_{x_i}\Big](u_h-\bu_h)(x)\leq 0.
\end{align*}
At the point $x$, we subtract the equations for $u_h$ and
$\bu_h$. After some rearranging using monotonicity of the scheme (the
above two inequalities) we get
\begin{align*}
0\leq
\sup_{\alp\in\A,\beta\in\B}\bigg\{&\sum_{i=1}^N\Big[(b_i^+-\bar
  b_i^+)(x) \Delta^+_{x_i} + (b_i^{-}-\bar b_i^-)(x) 
\Delta^-_{x_i}\Big]\bu_h(x)\\
&-c(x)(u_h-\bu_h)(x)-\bu_h(x)(c-\bar c)(x)+(f-\bar f)(x)\bigg\}.
\end{align*}
This (almost) immediately gives the upper bound on $u_h-\bu_h$.
Reversing the roles of $u_h$ and $\bu_h$ gives the lower bound and the
proof is complete.
\end{proof}


\begin{thebibliography}{10}

\bibitem{Am:OP} 
A.~L. Amadori. 
\newblock The obstacle problem for nonlinear integro-differential 
operators arising in option pricing. 
\newblock Quaderno IAC Q21-000, 2000.

\bibitem{BKR:BSAmOp}
F.~E. Benth, K.~H. Karlsen, and K. Reikvam.
A semilinear Black and Scholes partial differential equation for
valuing American options.
{\em Finance Stoch.} 7(3):277--298, 2003.

\bibitem{BJ:Rate}
G.~Barles and E.~R.~Jakobsen.
\newblock On the convergence rate of approximation schemes for
Hamilton-Jacobi-Bellman equations.
\newblock {\em M2AN Math. Model. Numer. Anal.} 36(1):33--54, 2002.

\bibitem{BJ:Rate2}
G.~Barles and E.~R.~Jakobsen.
Error bounds for monotone approximation schemes for
Hamilton-Jacobi-Bellman equations. Submitted.

\bibitem{BS:Conv} 
G.~Barles and P.~E.~Souganidis.
\newblock Convergence of approximation schemes for fully nonlinear
second order equations.
\newblock {\em Asymptotic Anal}.  4(3):271--283, 1991.

\bibitem{BL:Book}
A.~Bensoussan and J.-L.~Lions.
\newblock {\em Applications of Variational Inequalities in Stochastic
  Control.}
\newblock North-Holland Publishing Co., Amsterdam-New York, 1982. 

\bibitem{CF:Appr} 
F.~Camilli and M.~Falcone.
\newblock An approximation scheme for the optimal control of diffusion
processes.  
\newblock {\em RAIRO Mod{\'e}l. Math. Anal. Num{\'e}r.}  29(1),97--122.


\bibitem{CIL:UG}
M.~G. Crandall, H.~Ishii, and P.-L. Lions.
\newblock User's guide to viscosity solutions of second order partial
  differential equations.
\newblock {\em Bull. Amer. Math. Soc. (N.S.)}, 27(1):1--67, 1992.

\bibitem{CL:Appr}
M.~G. Crandall and P.-L. Lions.
\newblock Two approximations of solutions of Hamilton-Jacobi equations. 
\newblock {\em Math. Comp.} 43(167):1--19, 1984.

\bibitem{FS:Book}
W.~H. Fleming and H.~M. Soner.
\newblock {\em Controlled {M}arkov processes and viscosity solutions}.
\newblock Springer-Verlag, New York, 1993.

\bibitem{FS:SDG} W.~H. Fleming and P.~E. Souganidis. 
\newblock On the existence of value functions of two-player, zero-sum
   stochastic differential games. 
\newblock {\em Indiana Univ. Math. J.} 38(2):293--314, 1989.

\bibitem{Is:Equiv}
H.~Ishii.
\newblock On the equivalence of two notions of weak solutions,
viscosity solutions and distribution solutions.
\newblock {\em Funkcial. Ekvac.}, 38(1):101--120, 1995.

\bibitem{Is:Unique}
H.~Ishii.
\newblock On uniqueness and existence of viscosity solutions of fully nonlinear
  second-order elliptic {P}{D}{E}s.
\newblock {\em Comm. Pure Appl. Math.}, 42(1):15--45, 1989.

\bibitem{J:Non-conv}
E.~R.~Jakobsen.
\newblock Error bounds for monotone approximation schemes
  for non-convex degenerate elliptic equations in $\R^1$.
\newblock To appear in {\em BIT}.

\bibitem{J:Par}
E.~R. Jakobsen. 
\newblock On the rate of convergence of approximation schemes for
Bellman equations associated with optimal stopping time problems.
\newblock {\em Math. Models Methods Appl. Sci. (M3AS).} 13(5):
613-644, 2003. 


\bibitem{JK:ContDep}
E.~R. Jakobsen and K.~H. Karlsen.
\newblock Continuous dependence estimates for viscosity solutions of
fully nonlinear degenerate parabolic equations. 
\newblock {\em J. Differential Equations} 183:497-525, 2002.

\bibitem{JK:Ell}
E.~R. Jakobsen and K.~H. Karlsen.
\newblock Continuous dependence estimates for viscosity solutions of
fully nonlinear degenerate elliptic equations. 
\newblock {\em Electron. J. Diff. Eqns.} 2002(39):1--10, 2002. 

\bibitem{KS:Book}
D.~Kinderlehrer and G.~Stampacchia.
\newblock {\em An Introduction to Variational Inequalities and their
Applications.}
Reprint of the 1980 original. SIAM, Philadelphia, 2000. 

\bibitem{Kr:HJB1}
N.~V. Krylov.
\newblock On the rate of convergence of finite-difference approximations
for Bellman's equations.
\newblock {\em St. Petersburg Math. J.}, 9(3):639--650, 1997.

\bibitem{Kr:HJB2}
N.~V. Krylov.
\newblock On the rate of convergence of finite-difference approximations
for Bellman's equations with variable coefficients.
\newblock {\em Probab. Theory Ralat. Fields}, 117:1--16, 2000.


\bibitem{KD:Book}
H.~J. Kushner and P.~Dupuis.
\newblock {\em Numerical methods for for stochastic control
problems in continuous time.}
\newblock Springer-Verlag, New York, 2001.

\bibitem{Li:DPI1}
P.~L. Lions.
\newblock Optimal control of diffusion processes and Hamilton-Jacobi-Bellman
equations, Part I: The dynamic programming principle and applications.
\newblock Comm. P.D.E. {\bf 8} (1983).

\bibitem{Li:DPI2}
P.~L. Lions.
\newblock Optimal control of diffusion processes and Hamilton-Jacobi-Bellman
equations, Part II: Viscosity solutions and uniqueness.
\newblock Comm.P.D.E. {\bf 8} (1983).

\bibitem{Men}
J.~L. Menaldi
\newblock Some estimates for finite difference approximations.
\newblock {\em SIAM J. Control and Optimization} {\bf 27} No 3 (1989)
579-607.

\bibitem{Ph:OS}
H.~Pham.
\newblock Optimal stopping of controlled jump diffusion processes: A
viscosity solution approach.
\newblock  {\em J. Math. Systems Estim. Control} 8(1), 1998.

\bibitem{WHD:Book} 
P.~Wilmott, S.~Howison, and J.~Dewynne.
\newblock {\em The Mathematics of financial Derivatives. A student
    introduction.} 
\newblock Cambridge University Press, Cambridge, 1995.

\end{thebibliography}
\end{document}